\documentclass[a4paper,11pt]{article}
\usepackage{amsfonts}
\usepackage{bbm}
\usepackage{mathrsfs}
\usepackage{amsmath,amsthm,amssymb,amscd}
\usepackage[center]{titlesec}
\newtheorem{theo}{Theorem}[section]
\newtheorem{prop}[theo]{Proposition}
\newtheorem{lem}[theo]{Lemma}
\newtheorem{rem}[theo]{Remark}

\newtheorem{defi}[theo]{Definition}
\newcommand{\bp}{\begin{proof}}
\newcommand{\ep}{\end{proof}}

\begin{document}
\setlength{\baselineskip}{13pt} \pagestyle{myheadings}
\title{{ On a generalized Calabi-Yau equation}\thanks{ Supported by
 NSF Grant (China) 10671171.}}
\author{{\large Hongyu WANG\thanks {E-mail:
hywang@yzu.edu.cn} \ \ \ \ Peng ZHU }\thanks {E-mail:
zhupeng2004@yahoo.com.cn}\\
{\small School of Mathematical Science, Yangzhou University,}\\
{\small Yangzhou, Jiangsu 225002, P. R. China. }}
\date{}
\maketitle
\begin{quote}
\noindent {\bf Abstract.} {\small Dealing with the generalized
Calabi-Yau equation proposed by Gromov on closed almost-K\"{a}hler
manifolds, we extend to arbitrary dimension a non-existence result
proved in complex dimension $2$.}
\end{quote}

\section{Introduction}
The Calabi conjecture \cite{Ca} asserts that any representative of
the first Chern class of a closed K\"{a}hler manifold $(M,\omega)$
of real dimension $2n$ can be written as the Ricci curvature of a
K\"{a}hler metric $\omega'$ cohomologous to $\omega$. This
conjecture was proved by S. T. Yau \cite{Ya1}. Yau's result is
equivalent to finding a K\"{a}hler metric in a given K\"{a}hler
class with a prescribed volume form. More precisely, by
$\partial\bar{\partial}$-lemma on a K\"{a}hler manifold, this
amounts to solve the complex Monge-Amp\`{e}re equation:
\begin{align}\label{1000}
(\omega+\sqrt{-1}\partial\bar{\partial}\phi)^n=e^F\omega^n,
\end{align}
for some real function $\phi$ with
$$\omega+\sqrt{-1}\partial\bar{\partial}\phi>0,$$
 where $F\in C^\infty(M;\mathbb{R})$ satisfies that
\begin{align*}
\int_Me^F\omega^n=\int_M\omega^n.
\end{align*}
We call Equation \eqref{1000} Calabi-Yau equation. Yau solved this
equation by virtue of the continuity method. Consider the system of
equations obtained by
 replacing $F$ by $tF+c_t$, where $c_t$ is a constant for $t\in[0,1]$
 ($c_0=c_1=0$). An openness argument is obtained by the
 implicit function theorem and a closeness is derived by a priori estimate.
Similar questions are proposed in symplectic manifolds in different
cases and studied by many authors \cite{De,Do,ToWeYa,We}.\\

Suppose that $(M,\omega)$ is a closed $2n$-dimensional symplectic
manifold with a volume form $\sigma\in [\omega^n]$. In \cite{Mos},
Moser proved that there exists a symplectic form $\omega'\in
[\omega]$ which is symplectomorphic to $\omega$ satisfying that
\begin{align}\label{40001}
\omega'^n=\sigma.
\end{align}
That is, there exists a diffeomorphism,
$$f:M\rightarrow M,$$ isotopic to the identity such that
$$f^*\omega^n=\sigma.$$
It is easy to see that $\omega'=f^*\omega$ is cohomologous to
$\omega$ and $f^*\omega$ satisfies Equation \eqref{40001}. This is
independent of the almost complex structures. It is well known that
there are many almost complex structures which are compatible with
the symplectic form $\omega$ and they form a contractible space
\cite{Au,McSa}. Here $\omega$ is compatible with an almost complex
structure $J$, that is, at every point $p\in M$,  $\omega_p(v,Jv)>0$
for every nonzero vector $v\in T_pM$ and $\omega(JY,JZ)=\omega(Y,Z)$
for all vector
fields $Y$ and $Z$.\\

Now we suppose that $(M,g,J,\omega)$ is an almost K\"{a}hler
manifold, that is, the symplectic form $\omega$ is compatible with
the almost complex structure $J$ and $g(X,Y)=\omega(X,JY)$.
Obviously, $g$ is a Riemannian metric.
 Consider the existence of the solution of Equation
\eqref{40001} in the following form,
\begin{align}\label{1004}
\omega'=\omega(\phi)\equiv\omega+d(Jd\phi),
\end{align}
for $$\phi\in C^{\infty}(M;\mathbb{R}).$$ Here
$$(Jd\phi)(X)=d\phi(JX),$$ and  $\omega'$  tames $J$, that is,
at every point $p\in M$,  $\omega'_p(v,Jv)>0$ for every nonzero
vector
$v\in T_pM$. More precisely, this question can be expressed as follows:\\

 Does there exist a smooth function,
 $$\phi\in C^{\infty}(M;\mathbb{R}),$$  satisfying the
 following conditions?
\begin{equation}\label{40002}
\begin{cases}
\omega'^n=\sigma,\\
\omega'=\omega(\phi)\equiv\omega+dJd\phi\ {\rm tames} \ J,
\end{cases}
\end{equation}
where $\sigma$ is a given volume form in $[\omega^n]$. Following a
suggestion of M. Gromov, P. Delano\"{e} studied this problem in
\cite{De}.\\

 We call Equation \eqref{40002} generalized Calabi-Yau
equation.  In particular, if $(M,g,J,\omega)$ is K\"{a}hler (that
is, $J$ is integrable), then,
\begin{align*}
d(Jd\phi)&=d(J(\partial+\bar{\partial})\phi)\\
&=\sqrt{-1}(\partial+\bar{\partial})(\partial-\bar{\partial})\phi\\
&=2\sqrt{-1}\bar{\partial}\partial\phi.
\end{align*}
This implies  that the form $\omega'$ is a K\"{a}hler form. So
Equation \eqref{40002}  reduces to the Calabi-Yau equation on
K\"{a}hler manifolds, which was solved by S. T. Yau
\cite{Ya1}.\\

We define an operator $F$ from $C^\infty(M;\mathbb{R})$ to
$C^\infty(M;\mathbb{R})$ as follows:
\begin{align*}
\phi&\mapsto F(\phi),
\end{align*}
where
\begin{align}\label{1003}
F(\phi)\omega^n=(\omega(\phi))^n.
\end{align}
Therefore, Equation \eqref{40002} is equivalent to the following
problem:\\

Suppose that a positive function $f\in C^{\infty}(M;\mathbb{R})$
satisfies the following equality,
\begin{align*}
\int_M\omega^n=\int_Mf\omega^n.
\end{align*}
Does there exist a solution of $\phi\in C^{\infty}(M;\mathbb{R})$
which satisfies the following equation?
\begin{equation}\label{4000}
\begin{cases}
F(\phi)=f,\\
\omega(\phi)\equiv\omega+dJd\phi\ {\rm tames} \ J.
\end{cases}
\end{equation}\\

We need some notations in \cite{De}:
\begin{defi} Suppose that $(M,g,J,\omega)$ is an almost
K\"{a}hler manifold of real dimension $2n$. The sets $A$, $B$, $A_+$
and $B_+$ are defined as follows:
\begin{align*}
A&:=\{\phi\in C^{\infty}(M;\mathbb{R})\mid\int_M\phi\ \omega^n=0\};\\
 B&:=\{f\in C^{\infty}(M;\mathbb{R})\mid\int_Mf\omega^n=\int_M\omega^n\};\\
  A_+&:=A\cap\{\phi\in C^{\infty}(M;\mathbb{R})\mid\omega(\phi) \ {\rm tames} \ J\};\\
   B_+&:=B\cap\{f\in C^{\infty}(M;\mathbb{R})\mid f>0\}.
\end{align*}
\end{defi}
Note that $A_+$ can be regarded as a convex open set of symplectic
potential functions (analogue of K\"{a}hler potential functions).
Restricting the operator $F$ to $A_+$, we get
\begin{align*}
F(A_+)\subset B_+.
\end{align*}
Thus, the existence of a solution to Equation \eqref{4000} is
equivalent to that the restricted operator
\begin{align}\label{40004}
F|_{A_+}:A_+\rightarrow B_+,
\end{align}
is surjective.\\

Suppose that $(M,g,J,\omega)$ is a closed almost K\"{a}hler manifold
of real dimension $2n$. If $J$ is integrable, then $F:A_+\rightarrow
B_+$ is surjective.  Conversely, if $F:A_+\rightarrow B_+$  is
bijective and $n=2$, then $J$ is integrable. Delano\"{e} \cite{De}
proved this result by constructing a suitable smooth function
$\phi_0$ on the boundary of $A_+$ and he conjectured
\cite[Conjecture p.837]{De} that the same result holds when $n\geq2$.\\

In this paper, we prove Delano\"{e}'s conjecture. Since an oriented
surface is K\"{a}hler, we always consider the case $n\geq2$.
\begin{theo}\label{1001}
Suppose that $(M,g,J,\omega)$ is a closed almost K\"{a}hler manifold
of real dimension $2n$. Then the restricted operator
$F|_{A_+}:A_+\rightarrow F(A_+)$ is a diffeomorphism. Moreover, the
restricted operator $F|_{A_+}:A_+\rightarrow B_+$ is a surjectivity
map if and only if $J$ is integrable.
\end{theo}
\begin{rem}
{\rm 1)} S. K. Donaldson  gave a conjecture in \cite{Do}. Suppose
that $J$ is an almost complex structure on a closed symplectic
manifold $(M,\omega)$ of dimension $4$, and is tamed by the
symplectic form $\omega$. Let $\sigma$ be a smooth volume form on
$M$ with
$$\int_M\sigma=\int_M\omega^2.$$
 Let $\tilde{\omega}$ be an almost
K\"{a}hler form corresponding to the almost complex structure $J$
(that is, $\tilde{\omega}$ is a symplectic form compatible with $J$)
with $[\tilde{\omega}]=[\omega]$, and solve the Calabi-Yau equation
$$\tilde{\omega}^2=\sigma.$$
Donaldson conjectured that there are  $C^\infty$ a prior bounds on
$\tilde{\omega}$ depending only on $\omega$, $J$ and $\sigma$. This
is  related to his broader program
\cite{Do}.\\

{\rm 2)} V. Tosatti, B. Weinkove and S. Y. Yau \cite{ToWeYa} proved
that Donaldson's conjecture on estimates for Calabi-Yau equation in
terms of a taming symplectic form can be reduced to an integral
estimate of a scalar potential function. We will study Donaldson's
conjecture in a future paper \cite{WaZh}.
\end{rem}
This article is organized as follows. Section \ref{2} contains a
local theory of the operator $F$ (Equality \eqref{1003}) and the
local expression of $\omega(\phi)$ (Equality \eqref{1004}) by
choosing the second canonical connection. Section \ref{3} gives the
 proof of Theorem \ref{1001}. In the following, we
simply call Equation \eqref{40002} the  Calabi-Yau equation.

\section{Local theory of  Calabi-Yau equation}\setcounter{equation}{0}\label{2}
This section is devoted to a local theory of Calabi-Yau equation.
Let $(M,g,J,\omega)$ be a closed almost K\"{a}hler manifold of real
dimension $2n$. The $\mathbb{C}$-extension of the almost complex
structure $J$ on the complexified tangent bundle
$TM\otimes\mathbb{C}$ is defined by
$$J(X+iY)=JX+iJY$$
for $X,Y \in TM.$ Obviously,
$$J^2=-Id$$
 on $TM\otimes\mathbb{C}$.
Therefore,
$$TM\otimes\mathbb{C}=T^{1,0}\oplus T^{0,1},$$
where $$T^{1,0}=\{X\in TM\otimes\mathbb{C}:JX=iX\}.$$ $J$ induces an
almost complex structure on $\wedge^kT^*M\otimes\mathbb{C}$. Hence
$$\wedge^2T^*M\otimes\mathbb{C}=T_{2,0}+T_{1,1}+T_{0,2},$$
where $T_{p,q}$ is the space of $(p,q)$-forms. Let $P_{2,0}$,
$P_{1,1}$ and $P_{0,2}$
 be the projections to $T_{2,0}$, $T_{1,1}$ and $T_{0,2}$,
respectively. For each $\phi\in C^\infty(M;\mathbb{R})$, set
\begin{align}\label{20001}
\tau(\phi)=P_{2,0}(\omega(\phi)),
\end{align}
\begin{align}\label{20002}
H(\phi)=P_{1,1}(\omega(\phi)).
\end{align}
\begin{prop}\label{20000} Suppose that $(M,g,J,\omega)$ is a closed almost K\"{a}hler manifold of real dimension $2n$.
 Then $\omega(\phi)$ tames $J$ if and only if $H(\phi)$ tames $J$.
\end{prop}
\bp From the definition of $\omega(\phi)$ and $H(\phi)$, we have
\begin{align*}
\omega(\phi)(X,JX)&=(\tau(\phi)+H(\phi)+\overline{\tau(\phi)})(X,JX)\\
&=H(\phi)(X,JX).
\end{align*}
The last equality holds since $\tau(\phi)(X,JX)=0$, for every vector
field $X\in TM$.
 \ep
\begin{rem} Observe that the taming $(1,1)$-form $H(\phi)$ in
Proposition \ref{20000} is not necessarily closed; it is, when $J$
is integrable (K\"{a}hler case) in which case it coincides with
$\omega(\phi)$.
\end{rem}
To compute $\tau(\phi)$ and $H(\phi)$, we choose a local coordinate
system and the second canonical connection on an almost Hermitian
manifold, $(M,g,J,\omega)$ of real dimension $2n$, (that is, a
Riemannian metric $g$ is $J$-invariant and $\omega(X,Y)=g(JX,Y)$. ).

 There exists a local orthonormal basis $\{\epsilon_1,J\epsilon_1,\cdots,\epsilon_n,J\epsilon_n\}$
  of $M$. Set $e_j=\frac{\epsilon_j-\sqrt{-1}J\epsilon_j}{\sqrt{2}}$.
 Then $\{e_1,\cdots,e_n,\bar{e}_1,\cdots,\bar{e}_n\}$ is a local basis
  of $TM\otimes\mathbb{C}$ and $g(e_i,e_j)=g(\bar{e}_i,\bar{e}_j)=0,$
$g(e_i,\bar{e}_j)=\delta_{ij}.$ Obviously, $\{e_1,\cdots,e_n\}$
 is a local basis of $T^{1,0}.$ Let
$\{\theta^1,\cdots,\theta^n\}$ be its dual basis. Obviously, after
complexification, $J$, $g$, and $\omega$ can be expressed locally as
follows,
 \begin{align*}
J=\sqrt{-1}(\theta^{\alpha}\otimes
e_{\alpha}-\bar{\theta}^{\alpha}\otimes\bar{e}_{\alpha}),
\end{align*}
 \begin{align*}
g=\sum_{i=1}^n(\theta^{\alpha}\otimes
\bar{\theta}^{\alpha}+\bar{\theta}^{\alpha}\otimes\theta^{\alpha}),
\end{align*}
 \begin{align}\label{20010}
\omega=\sqrt{-1}\sum_{\alpha=1}^n\theta^{\alpha}\wedge\bar{\theta}^{\alpha}.
\end{align}
Now, although this choice is not essential for the proof, we choose
a different connection than the one used in \cite{De}, namely we
choose the second canonical connection $\nabla^1$ (an almost
Hermitian connection) on an almost K\"{a}hler manifold
\cite{EhLi,Ga,ToWeYa,Zh}. Note that the second canonical connection
on the almost Hermitian manifold $(M,g,J,\omega)$ is an affine
connection $\nabla^1$ satisfying
\begin{align}\label{2009}
\nabla^1 g=0=\nabla^1 J,
\end{align}
and the $(1,1)$-component of its torsion vanishes
\cite{Ga,ToWeYa,Zh}. We must redo with $\nabla^1$ calculations
analogous to those presented in \cite[Appendix 1]{De}. From Equality
\eqref{2009},  we obtain that
\begin{align*}
\nabla^1:T^{1,0}&\rightarrow (T^*M\otimes \mathbb{C})\otimes T^{1,0},\\
e_{\alpha}&\mapsto\omega_{\alpha}^{\beta}\otimes e_{\beta}.
\end{align*}
The second canonical connection, $\nabla^1$, induces
\begin{align*}
\nabla^1:T_{1,0}&\rightarrow (T^*M\otimes \mathbb{C})\otimes T_{1,0},\\
\theta^{\alpha}&\mapsto-\omega_{\beta}^{\alpha}\otimes
\theta^{\beta}.
\end{align*}
Let $\phi\in C^{\infty}(M;\mathbb{R})$. Then
\begin{align}\label{20011}
d\phi=\phi_{\alpha}\theta^{\alpha}+\bar{\phi}_{\alpha}\bar{\theta}^{\alpha}.
\end{align}
Thus,
\begin{align}\label{2010}
Jd\phi=\sqrt{-1}(\phi_{\alpha}\theta^{\alpha}-\bar{\phi}_{\alpha}\bar{\theta}^{\alpha}).
\end{align}
Let
\begin{align}\label{2011}
d\phi_{\alpha}-\phi_{\beta}\omega_{\alpha}^{\beta}
=\phi_{\alpha\beta}\theta^{\beta}+\phi_{\alpha\bar{\beta}}\bar{\theta}^{\beta}.
\end{align}
Then, equality $d^2\phi=0$ and Equation \eqref{20011} imply that
\begin{align}\label{2012}
d(\phi_{\alpha}\theta^{\alpha})+d(\bar{\phi}_{\alpha}\bar{\theta}^{\alpha})=0.
\end{align}

Let
\begin{align}\label{20012}
\Theta^\alpha=d\theta^\alpha+\omega_\beta^\alpha\wedge\theta^\beta
\end{align}
be the torsion of the second canonical connection $\nabla^1$. Thus,
$\Theta^\alpha$ contains $(2,0)$ and $(0,2)$ components only.
Therefore,
\begin{align}\label{20029}
\Theta^\alpha=T_{\beta\gamma}^\alpha\theta^\beta\wedge\theta^\gamma+
N_{\bar{\beta}\bar{\gamma}}^\alpha\bar{\theta}^\beta\wedge\bar{\theta}^\gamma,
\end{align}
with $T_{\beta\gamma}^\alpha=-T_{\gamma\beta}^\alpha$ and
$N_{\bar{\beta}\bar{\gamma}}^\alpha=-N_{\bar{\gamma}\bar{\beta}}^\alpha$.
Indeed, the $(0,2)$ component of the torsion is independent of the
choice of a metric and can be regarded as the Nijenhuis tensor, $N,$
of the almost complex structure $J$ \cite{ToWeYa,Zh}. From Equations
\eqref{20010}
 and \eqref{20012}, we obtain that
\begin{align}\label{20013}
d\omega&=d(\sqrt{-1}\sum_{\alpha=1}^n\theta^{\alpha}\wedge\bar{\theta}^{\alpha})\nonumber\\
&=\sqrt{-1}\sum_{\alpha=1}^n(d\theta^{\alpha}\wedge\bar{\theta}^{\alpha}
-\theta^{\alpha}\wedge d\bar{\theta}^{\alpha})\nonumber\\
&=\sqrt{-1}\sum_{\alpha=1}^n(\Theta^{\alpha}\wedge\bar{\theta}^{\alpha}
-\theta^{\alpha}\wedge \bar{\Theta}^{\alpha})\nonumber\\
&=\sqrt{-1}\sum_{\alpha,\beta,\gamma=1}^n
(T_{\beta\gamma}^\alpha\bar{\theta}^{\alpha}\wedge\theta^\beta\wedge\theta^\gamma
-\bar{T}_{\beta\gamma}^\alpha\theta^{\alpha}\wedge\bar{\theta}^\beta\wedge\bar{\theta}^\gamma\nonumber\\
&+N_{\bar{\beta}\bar{\gamma}}^\alpha\bar{\theta}^{\alpha}\wedge\bar{\theta}^\beta\wedge\bar{\theta}^\gamma
-\bar{N}_{\bar{\beta}\bar{\gamma}}^\alpha\theta^{\alpha}\wedge\theta^\beta\wedge\theta^\gamma).
\end{align}
Suppose that  $(M,g,J,\omega)$ is an almost K\"{a}hler manifold. By
Equation \eqref{20013}, we have
 $$T_{\beta\gamma}^\alpha=0,$$
 and
 $$N_{\bar{\alpha}\bar{\beta}}^\gamma+N_{\bar{\beta}\bar{\gamma}}^\alpha+N_{\bar{\gamma}\bar{\alpha}}^\beta=0.$$
Combining with  Equations \eqref{2010}$-$\eqref{20029}, we obtain
that
\begin{align*}
dJd\phi&=2\sqrt{-1}d(\phi_{\alpha}\theta^{\alpha})\nonumber\\
&=2\sqrt{-1}(d\phi_{\alpha}\wedge\theta^{\alpha}+\phi_{\alpha}d\theta^{\alpha})\nonumber\\
&=2\sqrt{-1}(\phi_{\alpha\beta}\theta^\beta\wedge\theta^\alpha
+\phi_{\alpha\bar{\beta}}\bar{\theta}^\beta\wedge\theta^\alpha+\phi_\alpha\Theta^\alpha)\nonumber\\
&=2\sqrt{-1}(\phi_{\alpha\beta}\theta^\beta\wedge\theta^\alpha
+\phi_{\alpha\bar{\beta}}\bar{\theta}^\beta\wedge\theta^\alpha
+\phi_\alpha
N_{\bar{\beta}\bar{\gamma}}^\alpha\bar{\theta}^\beta\wedge\bar{\theta}^\gamma).
\end{align*}

The reality of $dJd\phi$ implies that
\begin{align*}
\phi_{\alpha\bar{\beta}}=\overline{\phi_{\beta\bar{\alpha}}},
\end{align*}
and
\begin{align*}
\phi_{\alpha\beta}\theta^\beta\wedge\theta^\alpha
=-\overline{\phi_\alpha
N_{\bar{\beta}\bar{\gamma}}^\alpha\bar{\theta}^\beta\wedge\bar{\theta}^\gamma}.
\end{align*}
Set
\begin{align}\label{20017}
h(\phi)(X,Y)=\frac{1}{2}(H(\phi)(X,JY)+H(\phi)(Y,JX)).
\end{align}
Then,
\begin{align}\label{20018}
\tau(\phi)=-2\sqrt{-1}\bar{\phi}_\alpha
\bar{N}_{\bar{\beta}\bar{\gamma}}^\alpha\theta^\beta\wedge\theta^\gamma,
\end{align}
\begin{align}\label{20019}
H(\phi)=\sqrt{-1}(\delta_{\alpha\bar{\beta}}-2\phi_{\alpha\bar{\beta}})\theta^\alpha\wedge\bar{\theta}^\beta,
\end{align}
\begin{align}\label{20020}
h(\phi)=(\delta_{\alpha\bar{\beta}}-2\phi_{\alpha\bar{\beta}})(\theta^\alpha\otimes\bar{\theta}^\beta
+\bar{\theta}^\beta\otimes\theta^\alpha),
\end{align}
\begin{align*}
dJd\phi&=-2\sqrt{-1}(\bar{\phi}_\alpha
\bar{N}_{\bar{\beta}\bar{\gamma}}^\alpha\theta^\beta\wedge\theta^\gamma
+\phi_{\alpha\bar{\beta}}\theta^\alpha\wedge\bar{\theta}^\beta
-\phi_\alpha
N_{\bar{\beta}\bar{\gamma}}^\alpha\bar{\theta}^\beta\wedge\bar{\theta}^\gamma).
\end{align*}
Where $\delta_{\alpha\bar{\beta}}=1$ if $\alpha=\beta$; otherwise,
$\delta_{\alpha\bar{\beta}}=0$. Hence, we have the following lemma:
\begin{lem}\label{2005}Suppose that $(M,g,J,\omega)$ is a $2n$-dimensional almost K\"{a}hler
manifold, and $\phi$ is a smooth real function on $M$.
 Then,
\begin{align*}
\omega(\phi)&=\tau(\phi)+H(\phi)+\bar{\tau}(\phi),
\end{align*}
where $\tau(\phi)$ and $H(\phi)$ are given by Equations
\eqref{20018} and \eqref{20019}, respectively.
\end{lem}

Let $[\frac{n}{2}]$ denote the integral part of a positive number
$\frac{n}{2}$.
 Following \cite{De}, let us define operators $F_j$ ($j=0,1,\cdots,[\frac{n}{2}]$) as follows:
\begin{align*}
F_j: C^\infty(M;\mathbb{R})&\rightarrow C^\infty(M;\mathbb{R}),\\
\phi&\mapsto F_j(\phi),
\end{align*}
where
\begin{align}\label{2001}
 F_j(\phi)\omega^n=\frac{n!}{j!^2(n-2j)!}(\tau(\phi)\wedge\overline{\tau(\phi)})^j\wedge
 H(\phi)^{n-2j}.
\end{align}
Then
\begin{align*}
\sum_{j=0}^{[\frac{n}{2}]}F_j(\phi)\omega^n&=(\tau(\phi)+H(\phi)+\bar{\tau}(\phi))^n
 =F(\phi)\omega^n.
\end{align*}
The following result was proved in \cite[Proposition 5]{De}:
\begin{prop}\label{2008} Suppose that $(M,g,J,\omega)$ is a closed almost K\"{a}hler manifold of real dimension
$2n$. Then
$$F(\phi)=F_0(\phi)+\cdots+F_{[\frac{n}{2}]}(\phi).$$
 If $\phi\in A_+,$ then $F_0(\phi)>0$ and
 $F_j(\phi)\geq0$, for $j=1,2,\cdots,[\frac{n}{2}]$.\\
\end{prop}
For any $\phi\in A_+$, it is easy to see that the tangent space at
$\phi$, $T_{\phi}A_+$, is $A$. For
 $u\in T_{\phi}A_+=A$, define $L(\phi)(u)$ by,
$$L(\phi)(u)=\frac{d}{dt}F(\phi+tu)|_{t=0}.$$
For a local theory of Calabi-Yau equation, we should study the
tangent map of the restricted operator $F|_{A_+}$. The following
result was obtained in \cite[Proposition 1]{De}:
\begin{lem}\label{2006} Suppose that $(M,g,J,\omega)$ is a closed almost K\"{a}hler manifold
of dimension $2n$. Then the restricted operator
$$F:A_+\rightarrow B_+$$
is elliptic type on $A_+$. Moreover, the tangent map,
$$dF_{\phi}=L(\phi)$$
of $F$ at $\phi\in A_+$ is a linear elliptic differential operator
of second order without zero-th term.
\end{lem}
The following result was proved in \cite[Theorem 2]{De}:
\begin{prop}\label{2007}
Suppose that $(M,g,J,\omega)$ is a closed almost K\"{a}hler manifold
of dimension $2n$. Then the restricted operator
$$F:A_+\rightarrow F(A_+)$$
is a diffeomorphic map.
\end{prop}
\begin{rem}
From the definition of $F(\phi)$ (Equality \eqref{1003}), we see
that
$$F(0)\omega^n=(\omega(0))^n=\omega^n.$$
Thus, $F(0)=1$. Therefore, Proposition \ref{2007} implies that the
solution to Equation \eqref{40002} exists and is unique if $\sigma$
is a small perturbation of $\omega^n$.
\end{rem}
\section{Global theory of Calabi-Yau equation}\setcounter{equation}{0}\label{3}
In this section, we give a proof of the main theorem. Let
$(M,g,J,\omega)$ be a closed almost K\"{a}hler manifold of real
dimension $2n$. If the almost complex structure $J$ is integrable,
then the surjectivity of the restricted operator \eqref{40004} is
equivalent to the existence of solutions of Calabi-Yau equation on
K\"{a}hler manifolds which was solved by S. T. Yau \cite{Ya1}.
Suppose that $J$ is not integrable. We will prove that the
restricted operator \eqref{40004} is not surjective. For the case
$n=2$, it was proved by P. Delano\"{e} \cite{De}. Now we give the
proof for the case $n\geq2$. More precisely, as in \cite[p.835]{De},
we want to construct a function $\phi_0$ belonging to the boundary
of the convex open set $A_+$ such that $F(\phi_0)\in B_+$. First,
let us give the definition of a pseudo holomorphic function
\cite{Au1,Au,Gr}.
\begin{defi}
A smooth complex-valued function $f$ on a manifold $M$ is a pseudo
holomorphic function at some point $p\in M$ if $df\circ J=idf$ at
$p$.
\end{defi}
To construct $\phi_0$ belonging to the boundary of $A_+$ and
satisfying inequality $F(\phi_0)>0$, we need the following lemmas:
\begin{lem}\label{3001}
 Let $(M,J)$ be an almost complex manifold. Suppose that $f$ is
a pseudo holomorphic function at some point $p\in M$. Then, for all
vector fields $X,Y$,
\begin{align*}
df(N(X,Y))=0,
\end{align*}
at $p$, where $N$ is the Nijenhuis tensor.
\end{lem}
\bp If $X,Y\in T^{1,0}$, then
\begin{align*}
N(X,Y)&=[JX,JY]-J[X,JY]-J[JX,Y]-[X,Y]\\
&=-2[X,Y]-2iJ[X,Y].
\end{align*}
Thus
\begin{align}\label{30000}
JN(X,Y)&=-2J[X,Y]-2iJ^2[X,Y]\nonumber\\
&=-iN(X,Y).
\end{align}
If $$df\circ J=idf,$$ then
$$df\circ J(N(X,Y))=idf(N(X,Y)).$$
Combining with Equality \eqref{30000}, we get
$$df(N(X,Y))=0,$$
at $p\in M$. If $X,Y\in T^{0,1}$, then
$$df(N(X,Y))=\overline{df(N(\bar{X},\bar{Y}))}=0,$$
at $p\in M$. If $X\in T^{0,1}$ and  $Y\in T^{1,0}$, then $N(X,Y)=0.$
Hence,
$$df(N(X,Y))=0,$$
at $p\in M$. The proof is completed.
 \ep
\begin{rem}
Note that if there exist $n$ pseudo holomorphic functions on a real
$2n$-dimensional almost Hermitian manifold $(M,g,J)$ which are
independent at some point $p\in M$, then the Nijenhuis tensor $N$
identically vanishes at $p$. This means that an integrable complex
structure is one with many (pseudo)-holomorphic functions. It is a
hard theorem (Newlander-Nirenberg integrability theorem for almost
complex structures) that the converse is also true. In general, an
almost complex manifold has no holomorphic functions at all. On the
other hand, it has a lot of pseudo-holomorphic curves (i.e., maps
$f:\mathbb{C}\rightarrow (M^{2n},g,J)$ such that $df\circ i=J\circ
df$) \cite{Gr}.
\end{rem}

Let $(M,g,J,\omega)$  be a closed almost K\"{a}hler manifold of real
dimension $2n$, where the almost complex structure $J$ is not
integrable. Hence there exists a neighborhood $U_0$ of some point
$p_0\in M$ where the Nijenhuis tensor $N$ does not vanish. Pick
$f=\phi_1+i\phi_2$, with $\phi_1$ and $\phi_2$ real-valued, a
complex function on $M$ which is not holomorphic on $U_0$. By Lemma
\ref{3001}, it satisfies
 \begin{align*}
df(N(X,Y))\neq0,
\end{align*}
on the neighborhood $U_0$. Without loss generality, we can suppose
that for some $1\leq\alpha<\beta \leq n$,
 \begin{align*}
d\phi_1(N(e_\alpha,e_\beta))=d\phi_1(e_{\bar{\gamma}}
N_{\alpha\beta}^{\bar{\gamma}})\neq0,
\end{align*}
on the neighborhood $U_0$ of $p_0$. By Equality \eqref{20020}, there
exists a suitable (non-zero) constant $c$ such that $h(c\phi_1)$,
defined in Section \ref{2} (Equalities  \eqref{20017} and
\eqref{20020}), is a positive definite Hermitian matrix on $M$. Set
$\phi=c\phi_1$. Then $\tau(\phi)$, defined in Section \ref{2}
(Equalities \eqref{20001} and \eqref{20018}) is non-zero on the
neighborhood $U_0$ of $p_0$. So we obtain the following lemma:
\begin{lem}\label{3002}Suppose that $(M,g,J,\omega)$ is a closed almost K\"{a}hler
manifold of real dimension $2n$, where the almost complex structure
$J$ is not integrable. There exists a function $\phi\in
C^\infty(M;\mathbb{R})$
 satisfying the following conditions:\\
{\rm (1)} $h(\phi)$ is a positive definite Hermitian matrix on $M$;\\
{\rm (2)} there is a neighborhood $U_0$ of a point $p_0\in M$ such
that $\tau(\phi)$ is non-zero on the neighborhood $U_0$.
 \end{lem}
 We can suppose that
 \begin{align}\label{30001}
H(\phi)=\sqrt{-1}\lambda_\alpha(0)\theta^\alpha\wedge\bar{\theta}^\alpha
\end{align}
 at $p_0$ \cite{Wey}. By Lemma \ref{3002}, we have that
$\lambda_\alpha(0)>0$, for $\alpha=1,\cdots,n$.
There exists a local coordinate system $\varphi^{-1}:U_0\rightarrow
\mathbb{C}^n$ such that $\varphi(0)=p_0$ and
 \begin{align}\label{30002}
\frac{\partial}{\partial z_\alpha}=e_\alpha
\end{align}
at $p_0$, for $\alpha=1,\cdots,n$.  Without loss generality, we can
suppose that there is a polydisk $\triangle$
  with center $p_0$:\\
$$\triangle=\{(z_1,\cdots,z_n)\in \mathbb C^n:|z_i|\leq1, 1\leq i \leq n\}\subset \varphi^{-1}(U_0)$$
and there is $\epsilon_1>0$ such that, on $\varphi(\triangle)$,
\begin{align}\label{30010}
|\tau_{12}(\phi)|\geq\epsilon_1>0.
\end{align}

Let $\tilde{\eta}(r)$ be a $C^2$ cut-off function such that
\begin{equation}\label{302}
\tilde{\eta}(r)=
\begin{cases}
1 \ \ \ r\leq\frac{1}{2},\\
0 \ \ \ r\geq1.
\end{cases}
\end{equation}
Given a real number $R>2$, let
\begin{align}\label{30012}
\Phi_R(z)=\sum_{i=1}^2\frac{\lambda_i(0)}{2}|z_i|^2\cdot\tilde{\eta}(|R^2z_i|)
\end{align}
and
\begin{align}\label{30014}
\eta_R(z)=\tilde{\eta}(|Rz_1|)\cdots\tilde{\eta}(|Rz_n|).
\end{align}
Define a function $\psi_{R}$ as follows,
\begin{align}\label{30009}
\psi_R(z)=\Phi_R(z)\eta_R(z),
\end{align}
where $z\in\triangle$. In the sequel, for short, we will abusively
denote by $z_i$, $\eta_R$, $\Phi_R$, $\psi_R$, $\sigma$, the
corresponding functions on $U_0\subset M$ obtained by composition
with $\varphi^{-1}$. With this abuse of notation, recalling
\eqref{2011}, we obtain by a direct calculation that
\begin{align}\label{303}
P_{1,1}(dJd(\psi_R))&=-2\sqrt{-1}\sum_{1\leq\alpha,\beta\leq n}
(\psi_R)_{\alpha\bar{\beta}}\theta^\alpha\wedge\bar{\theta}^\beta,
\end{align}
with
\begin{align}\label{30023}
(\psi_R)_{\alpha\bar{\beta}}=\bar{e}_\beta
e_\alpha\psi_R-(\psi_R)_\gamma\Gamma_{\alpha\bar{\beta}}^\gamma,
\end{align}
where
$\Gamma_{\alpha\bar{\beta}}^\gamma=\omega_\alpha^\gamma({\bar{e}}_\beta)$
are connection coefficients of the second canonical connection
bounded on $\triangle$. Set
$$\triangle_R=\{(z_1,\cdots,z_n):|z_i|\leq\frac{1}{R}, 1\leq i \leq n\}.$$
Obviously, $\triangle_R\subset\triangle$. Henceforth, we will freely
denote by the same letter $C$ various positive constants independent
of $R>2$.
\begin{lem}\label{3003}Suppose that $(M,g,J,\omega)$ is a closed almost K\"{a}hler manifold
of real dimension $2n$, and $\tilde{\eta}$, $\psi_R$,
 $\triangle$ and $U_0$ are
defined as before. There exists a constant $C$ (depending only on
$\lambda_1(0)$, $\lambda_2(0)$, $\tilde{\eta}$, $\omega$, $g$ and
$J$)
 such that,\\
{\rm (1)} on $\triangle$, for $1\leq\alpha\leq n$,
\begin{align*}
|(\psi_R)_{\alpha}|\leq\frac{C}{R^2};
\end{align*}
{\rm (2)} on $\triangle$, for $1\leq\alpha=\beta\leq 2$,
\begin{align*}
|(\psi_R)_{\alpha\bar{\beta}}|\leq C;
\end{align*}
otherwise,
\begin{align*}
|(\psi_R)_{\alpha\bar{\beta}}|\leq C(\sigma+\frac{1}{R}),
\end{align*}
where $\sigma$ is a continuous function (independent of $R$) and
$\sigma(0)=0$. Moreover, at the center of $\triangle$,
$|(\psi_R)_{\alpha\bar{\beta}}|=\frac{\lambda_\alpha(0)}{2}$, for
$1\leq\alpha=\beta\leq 2$; otherwise,
$|(\psi_R)_{\alpha\bar{\beta}}|=0$.
 \end{lem}
\bp From \eqref{302},
 we can see that $\psi_R$ is zero on $\triangle\setminus\triangle_R$.
 Hence it suffices to prove that the desired estimate holds on
 $\triangle_R$. Note that, for $1\leq\alpha,\beta,i\leq n$,
 the functions $|e_\alpha z_i|$, $|e_\alpha {\bar{z}}_i|$,
 $|{\bar{e}}_\beta z_i|$, $|{\bar{e}}_\beta {\bar{z}}_i|$,
 $|{\bar{e}}_\beta e_\alpha z_i|$ and  $|{\bar{e}}_\beta e_\alpha {\bar{z}}_i|$
 are continuous and bounded on $\Delta$.
From \eqref{30012} and \eqref{30014}, we have on $\triangle_R$
\begin{align}\label{30016}
|\eta_R|\leq 1, \ \ \ |\Phi_R|\leq C R^{-4}.
\end{align}
A direct computation shows that
\begin{align*}
e_\alpha\eta_R(z)=\sum_{i=1}^n&\big\{\tilde{\eta}(|Rz_1|)\cdots
\widehat{\tilde{\eta}(|Rz_i|)}\cdots\tilde{\eta}(|Rz_n|)\nonumber\\
& \cdot\tilde{\eta}'(|Rz_i|)\cdot[R(2|z_i|)^{-1}((e_\alpha
z_i)\bar{z}_i+(e_\alpha \bar{z}_i)z_i) ]\big\},
\end{align*}
where the hat means that the term is omitted. It is easy to see that
\begin{align}\label{30020}
|e_\alpha\eta_R|\leq CR,\ \ \ \ |{\bar{e}}_\alpha\eta_R|\leq CR.
\end{align}
Since
\begin{align*}
e_\alpha\Phi_R(z)=4^{-1}\sum_{i=1}^2&\big\{\lambda_i(0)[(e_\alpha
z_i)\bar{z}_i+(e_\alpha {\bar{z}}_i)z_i]\nonumber\\
 &\cdot[2\tilde{\eta}(|R^2z_i|) +R^2|z_i|\tilde{\eta}'(|R^2z_i|)]\big\},
\end{align*}
we obtain that
\begin{align}\label{30021}
|e_\alpha\Phi_R|\leq CR^{-2},\ \ \ |{\bar{e}}_\alpha\Phi_R|\leq
CR^{-2}.
\end{align}
Thus
\begin{align}\label{30022}
|(\psi_R)_\alpha|&\leq|(e_\alpha\Phi_R)\eta_R|+|\Phi_R(e_\alpha\eta_R)|\leq
CR^{-2}.
\end{align}
 Now we consider
\begin{align}\label{30024}
\bar{e}_\beta e_\alpha\psi_R= &\Phi_R(\bar{e}_\beta
e_\alpha\eta_R)+(e_\alpha\Phi_R)(\bar{e}_\beta\eta_R)\nonumber\\
&+(\bar{e}_\beta\Phi_R)(e_\alpha\eta_R)+(\bar{e}_\beta
e_\alpha\Phi_R)\eta_R.
\end{align}
A direct computation shows that
\begin{align*}
\bar{e}_\beta e_\alpha\eta_R(z)
 &=\sum_{1\leq i\neq j\leq
n}\big\{\tilde{\eta}(|Rz_1|)\cdots\widehat{\tilde{\eta}(|Rz_i|)}
\cdots\widehat{\tilde{\eta}(|Rz_j|)}\cdots\tilde{\eta}(|Rz_n|)\nonumber\\
&\ \ \ \ \ \ \ \ \ \
\cdot\tilde{\eta}'(|Rz_i|)\cdot[R(2|z_i|)^{-1}((e_\alpha
z_i)\bar{z}_i+(e_\alpha\bar{z}_i)z_i)]\nonumber\\
&\ \ \ \ \ \ \ \ \ \
\cdot\tilde{\eta}'(|Rz_j|)\cdot[R(2|z_j|)^{-1}(\bar{e}_\beta
z_j)\bar{z}_j+(\bar{e}_\beta \bar{z}_j)z_j)]\big\}\nonumber\\
&+\sum_{1\leq i\leq
n}\big\{\tilde{\eta}(|Rz_1|)\cdots\widehat{\tilde{\eta}(|Rz_i|)}\cdots\tilde{\eta}(|Rz_n|)\nonumber\\
&\ \ \ \ \ \
\cdot\big[\tilde{\eta}''(|Rz_i|)R^2(2|z_i|)^{-2}\cdot((e_\alpha
z_i)\bar{z}_i+(e_\alpha\bar{z}_i)z_i)\nonumber\\
&\ \ \ \ \ \ \ \ \cdot((\bar{e}_\beta
z_i)\bar{z}_i+(\bar{e}_\beta\bar{z}_i)z_i)+4^{-1}R\tilde{\eta}'(|Rz_i|)\nonumber\\
&\ \ \ \ \ \ \ \ \cdot \big(2(\bar{e}_\beta e_\alpha
z_i)|z_i|^{-1}\bar{z}_i+(e_\alpha z_i)(\bar{e}_\beta\bar{z}_i|z_i|^{-1}-(\bar{e}_\beta z_i)|z_i|^{-3}{\bar{z}_i}^2)\nonumber\\
&\ \ \ \ \ \ \ \ \ +2(e_{\bar{\beta}}e_\alpha
\bar{z}_i)|z_i|^{-1}z_i+(e_\alpha\bar{z}_i)(\bar{e}_\beta
z_i|z_i|^{-1}-(\bar{e}_\beta\bar{z}_i)|z_i|^{-3}{z_i}^2)\big)\big]
\big\}.
\end{align*}
Thus,
\begin{align}\label{30033}
|\bar{e}_\beta e_\alpha\eta_R|\leq CR^2.
\end{align}
A direct computation shows that
\begin{align}\label{30030}
\bar{e}_\beta e_\alpha\Phi_R
 =4^{-1}\sum_{i=1}^2&\lambda_i(0)\big\{\big[(\bar{e}_\beta e_\alpha
 z_i)\bar{z}_i
 +(e_\alpha z_i)(\bar{e}_\beta\bar{z}_i)
 +(\bar{e}_\beta e_\alpha \bar{z}_i)z_i
 \nonumber\\
 &+(e_\alpha \bar{z}_i)(\bar{e}_\beta z_i)\big]\cdot\big[2\tilde{\eta}(|R^2z_i|)
+R^2|z_i|\tilde{\eta}'(|R^2z_i|)\big] \nonumber\\
& +2^{-1}\big[(e_\alpha z_i)\bar{z}_i+(e_\alpha
{\bar{z}}_i)z_i\big]\cdot\big[(\bar{e}_\beta
z_i)\bar{z}_i+(\bar{e}_\beta\bar{z}_i)z_i\big]\nonumber\\
&\cdot\big[3R^2|z_i|^{-1}\tilde{\eta}'(|R^2z_i|)+R^4\tilde{\eta}''(|R^2z_i|)\big]\big\}.
\end{align}
If $\alpha=\beta=1$ or $\alpha=\beta=2$, then
\begin{align}\label{30031}
|\bar{e}_\beta e_\alpha\Phi_R|\leq C;
\end{align}
If $3\leq\alpha=\beta\leq n$ or $1\leq\alpha\neq\beta\leq n$, then
\begin{align}\label{30032}
|\bar{e}_\beta e_\alpha\Phi_R|\leq C(\sigma+R^{-2})\leq
C(\sigma+R^{-1}),
\end{align}
where $\sigma$ is the continuous function on $\triangle$ such that
 $\sigma(0)=0$, defined by:
\begin{align*}
\sigma(z)=\max\{|(e_\alpha z_i)(\bar{e}_\beta z_i)|(z),|(e_\alpha
z_i)(\bar{e}_\beta\bar{z}_i)|(z),|(e_\alpha \bar{z}_i)(\bar{e}_\beta
z_i)|(z),\\
|(e_\alpha \bar{z}_i)(\bar{e}_\beta \bar{z}_i)|(z): i=1,2,
 1\leq\alpha\neq\beta\leq n\ \
\mathrm{or}\ \  3\leq\alpha=\beta\leq n \}.
\end{align*}
 From
\eqref{30016}-\eqref{30032}, we obtain that
\begin{align*}
|(\psi_R)_{\alpha\bar{\beta}}|\leq \left\{
\begin{array}{ll}
C, & 1\leq\alpha=\beta\leq 2;\\
C(\sigma+\frac{1}{R}), & \mathrm{otherwise},
\end{array}
\right.
\end{align*}
where $\sigma$ is a continuous function and $\sigma(0)=0$. Besides,
by \eqref{30002} combined with \eqref{302} and \eqref{30030}, we
obtain at $z=0$ that
$|(\psi_R)_{\alpha\bar{\alpha}}|(0)=\frac{\lambda_\alpha(0)}{2}$,
for $1\leq\alpha\leq 2$; otherwise,
$|(\psi_R)_{\alpha\bar{\beta}}|(0)=0$. \ep
\begin{lem}\label{3004}Suppose that $(M,g,J,\omega)$ is a closed
almost K\"{a}hler manifold
 of real dimension $2n$, where the almost complex structure $J$ is not
 integrable. Suppose
 $\phi$, $\psi_{R}$,  $\varphi$, $\tilde{\eta}$, $\triangle$ and $U_0$ are defined as before.
 Let $s\in [0,1]$.  Then
 there exists $R_1>2$ independent of $s\in [0,1]$ such that,
 for each $R\geq R_1$, on $\varphi(\triangle_R)$,
\begin{align*}
|\tau_{12}(\phi+s\psi_{R}\circ\varphi^{-1})|\geq\frac{\epsilon_1}{2}>0,
\end{align*}
where $\epsilon_1>0$ is the constant occuring in \eqref{30010}.
\end{lem}
\bp Combining Lemma \ref{3003} with
$$\tau_{\alpha\beta}(\psi_{R}\circ\varphi^{-1})
=-2\sqrt{-1}({\psi}_{R}\circ\varphi^{-1})_{\gamma}\bar{N}_{\bar{\alpha}\bar{\beta}}^\gamma,$$
on $\varphi(\triangle_R)$,
 we have the following inequality,
$$|\tau_{\alpha\beta}(\psi_{R}\circ\varphi^{-1})|\leq\frac{C}{R^2}\leq\frac{C}{R}.$$
By Lemma \ref{3002} and Inequality \eqref{30010}, we obtain,
\begin{align*}
|\tau_{12}(\phi+s\psi_{R}\circ\varphi^{-1})|&\geq|\tau_{12}(\phi)|
-s|\tau_{12}(\psi_{R}\circ\varphi^{-1})|\\
&\geq|\tau_{12}(\phi)|
-|\tau_{12}(\psi_{R}\circ\varphi^{-1})|\\
&\geq\epsilon_1-\frac{C}{R}.
\end{align*}
We can choose $R_1> 2$ such that
$\epsilon_1-\frac{C}{R_1}\geq\frac{\epsilon_1}{2}$. Thus, for each
$R\geq R_1$, the desired estimate follows. \ep

For $s\in [0,1]$, we consider the function
$\phi+s\psi_{R}\circ\varphi^{-1}$. If $s=0$, then
\begin{align}\label{30028}
h(\phi+0\cdot\psi_{R}\circ\varphi^{-1})=h(\phi)
\end{align}
 is a positive
definite Hermitian metric on $M$. As observed above, we have
\begin{align}\label{30027}
P_{1,1}(dJd(\psi_{R}\circ\varphi^{-1}))(p_0)=-\sqrt{-1}\sum_{1\leq\alpha\leq2}
\lambda_{\alpha}(0)\theta^\alpha\wedge \bar{\theta}^\alpha,
\end{align}
which, recalling \eqref{30001}, yields:
\begin{align}\label{30029}
h(\phi+1\cdot\psi_{R}\circ\varphi^{-1})\circ\varphi(0)
=\sum_{3\leq\alpha\leq n}\lambda_\alpha(0)(\theta^\alpha\otimes
\bar{\theta}^\alpha +\bar{\theta}^\alpha\otimes \theta^\alpha).
\end{align}
The latter is a semi-positive definite Hermitian matrix at the
center of $\triangle_R$ (i.e., $p_0=\varphi(0)\in M$). Following
\cite[Definition 3]{De}, let us recall the notion of positivity
amplitude:
\begin{defi}
Suppose that $(M,g,J,\omega)$ is a closed almost Hermitian manifold
of real dimension $2n$. For each non-constant function $\phi\in
C^\infty(M;\mathbb{R})$, the positivity amplitude of $\phi$ is the
real number
\begin{align*}
a_{\omega}(\phi)=\sup\{&s\in(0,\infty):\ h(s\phi)
 {\rm \ is\  a \
positive  \ Hermitian \  matrix}\\
 &{\rm
 function\ on\ M }\},
\end{align*}
where $h(\phi)$ is defined in Section \ref{2} (Equalities
\eqref{20017} and \eqref{20020}).
\end{defi}
  Since $M$ is closed and the support of $\psi_{R}$ is compact,
the positivity amplitude is finite. Hence, by Lemmas
\ref{3002}$-$\ref{3004} and Equalities \eqref{30028}, \eqref{30029},
we obtain the following proposition:
\begin{prop}\label{3005}Suppose that $(M,g,J,\omega)$ is a closed almost K\"{a}hler manifold
of real dimension $2n$, where $J$ is not integrable. Suppose
$\triangle_R$, $\psi_{R}$, $R_1$, $\phi$, $\varphi$,
 $U_0$ ,$\triangle$
 and $a_{\omega}(\phi)$ are defined as before and set
$$\phi_R=\phi+a_{\omega(\phi)}(\psi_{R}\circ\varphi^{-1})\cdot\psi_{R}\circ\varphi^{-1}.$$
Then, for each $R\geq
R_1$, we have:\\
(1) $0< a_{\omega(\phi)}(\psi_{R}\circ\varphi^{-1})\leq1$;\\
(2) $\phi_R$ lies on the boundary of $A_+$;\\
(3)
$h(\phi+a_{\omega(\phi)}(\psi_{R}\circ\varphi^{-1})\cdot\psi_{R}\circ\varphi^{-1})$
is a semi-positive definite Hermitian matrix function on
$\varphi(\triangle_R)$ and a positive definite Hermitian matrix
function on $M\setminus\varphi(\triangle_R)$.
 \end{prop}
Now we will prove the following key lemma, sticking to the
assumptions of Proposition \ref{3005}:
\begin{lem}\label{3006}
There exists $R_0\geq R_1$ large enough such that the function
$\phi_0=\phi_{R_0}$ satisfies $F(\phi_0)>0$ on $M$.
\end{lem}
\bp
 Note that
\begin{align}\label{30025}
H(\phi_R)=H(\phi)+a_{\omega(\phi)}(\psi_{R}\circ\varphi^{-1})\cdot
P_{1,1}\big(dJd(\psi_{R}\circ\varphi^{-1})\big).
\end{align}
From  \eqref{2001}, we have at $p_0$
\begin{align*}
F_1(\phi_R)=\sum_{1\leq k\neq l\leq n}|\tau_{kl}(\phi_R)|^2\cdot
 \det\big((h(\phi_R)_{\alpha\bar{\beta}})_{\alpha,\beta\in\{1,\ldots,\hat{k},\ldots,\hat{l},\ldots,n\}}\big).
 \end{align*}
From Equalities \eqref{30001}, \eqref{30027}, Lemma \ref{3004} and
Proposition \ref{3005}, we have,  at $p_0$, that
\begin{align*}
F_1(\phi_R)&\geq |\tau_{12}(\phi_R)|^2\cdot
 \det((h(\phi_R)_{\alpha\bar{\beta}})_{3\leq\alpha,\beta\leq
 n})\geq\big(\frac{\epsilon_1}{2}\big)^2\prod_{3\leq\alpha\leq n}\lambda_\alpha(0)>0.
\end{align*}
Let
\begin{align*}
H(\phi_R)=H'+H'',
\end{align*}
where
\begin{align*}
H'&=\sqrt{-1}\sum_{1\leq\alpha\neq\beta\leq
n}H(\phi)_{\alpha\bar{\beta}}\theta^\alpha\wedge
\bar{\theta}^\beta+\sqrt{-1}a_{\omega(\phi)}(\psi_{R}\circ\varphi^{-1})\\
&\cdot\big[\sum_{1\leq\alpha\neq\beta\leq
n}(\psi_R)_{\alpha\bar{\beta}}\theta^\alpha\wedge
\bar{\theta}^\beta+\sum_{3\leq\alpha\leq
n,}(\psi_R)_{\alpha\bar{\alpha}}\theta^\alpha\wedge
\bar{\theta}^\alpha\big]
\end{align*}
and
\begin{align*}
H''=\sqrt{-1}\sum_{1\leq\alpha\leq2}H(\phi_R)_{\alpha\bar{\alpha}}\theta^\alpha\wedge
\bar{\theta}^\alpha+ \sqrt{-1}\sum_{3\leq\alpha\leq
n}H(\phi)_{\alpha\bar{\alpha}}\theta^\alpha\wedge
\bar{\theta}^\alpha.
\end{align*}
Note that $H(\phi)_{\alpha\bar{\alpha}}(p_0)=\lambda_\alpha(0)$, for
$1\leq\alpha\leq n$. $H(\phi)_{\alpha\bar{\beta}}(p_0)=0$, for
$\alpha\neq\beta$. Let us define a neighborhood of $0\in \triangle$
by:
\begin{align*}
V_\delta=\{z\in
\triangle:|\sigma|(z)<\delta,|H(\phi)_{\alpha\bar{\alpha}}|(\varphi(z))>\frac{\lambda_\alpha(0)}{2},\\
 |H(\phi)_{\beta\bar{\gamma}}|(\varphi(z))<\delta,
 \beta\neq\gamma\},
\end{align*}
for some $\delta\in (0,\frac{1}{2})$. Recalling \eqref{2001}, we
have
\begin{align*}
F_1(\phi_R)&=n(n-1)\tau(\phi_R)\wedge\overline{\tau(\phi_R)}\wedge(H'+H'')^{n-2}/\omega^n\\
&=n(n-1)\{\tau(\phi_R)\wedge\overline{\tau(\phi_R)}\wedge(H'')^{n-2}\\
&+\tau(\phi_R)\wedge\overline{\tau(\phi_R)}\wedge\sum_{k=1}^{n-2}C_{n-2}^k(H')^k\wedge(H'')^{n-2-k}\}/\omega^n.
\end{align*}
 By Lemma \ref{3003}, we obtain that, on $\triangle_R\bigcap V_\delta$, each component of
 $H'$ is less than $C(\delta+\frac{1}{R})$
and $|(\tau(\phi_R))_{\alpha\beta}|\leq C$, where $C$ is independent
of $R$. By Lemma \ref{3002}, we have that, on $\triangle_R\cap
V_\delta$,
\begin{align*}
0<H(\phi)_{\alpha\bar{\alpha}}\leq C,
\end{align*}
where $C$ is independent of $R$ and $\alpha=3,\cdots,n$. By Lemma
\ref{3003} and  Proposition \ref{3005}, we have that, on
$\triangle_R\cap V_\delta$,

\begin{align*}
0\leq H(\phi_R)_{\alpha\bar{\alpha}}&\leq
H(\phi)_{\alpha\bar{\alpha}}+|(\psi_R)_{\alpha\bar{\alpha}}|\leq C ,
\end{align*}
where $\alpha=1,2$ and $C$ is independent of $R$. Therefore, on
$\triangle_R\cap V_\delta$, recalling the positivity lemma of
\cite[Appendix 3]{De}, we have
\begin{align*}
F_1(\phi_R)&\geq
|(\tau(\phi_R))_{12}|^2\Pi_{\alpha=3}^n(H(\phi)_{\alpha\bar{\alpha}})
-C(\delta+\frac{1}{R})\\
&\geq2^{-n}\epsilon_1^2\Pi_{\alpha=3}^n\lambda_\alpha(0)-C(\delta+\frac{1}{R})\\
&\equiv\epsilon_0-C(\delta+\frac{1}{R}),
\end{align*}
where $C$ and $\epsilon_0$ is independent of $R$. Choose
\begin{align*}
\delta_1<\min\{\frac{\epsilon_0}{4C},\frac{1}{2}\},\ \ \ \
R_2>\max\{R_1,\frac{4C}{\epsilon_0}\}.
\end{align*}
 Thus, for each $R>R_2$,
\begin{align*}
F_1(\phi_R)\geq\frac{\epsilon_0}{2}>0,
\end{align*}
 on $V_{\delta_1}\bigcap\triangle_R$. Moreover, we can choose $R_0>R_2$
such that $\triangle_{R_0}\subset V_{\delta_1}$. Define
\begin{align}\label{30015}
\phi_0=\phi+a_{\omega(\phi)}(\psi_{R_0}\circ\varphi^{-1})\cdot\psi_{R_0}\circ\varphi^{-1}.
\end{align}
By construction, we have $\phi_0=\phi$ on
$M\setminus\varphi(\triangle_{R_0})$ and, on
$\varphi(\triangle_{R_0})$:
\begin{align*}
F_1(\phi_0)\geq\frac{\epsilon_0}{2}>0,\ \ \ \ F_0(\phi_0)\geq0,
\end{align*}
and also \cite[Appendix 3]{De}: $\forall j\geq2, F_j(\phi_0)\geq0$.
Therefore, by  Proposition \ref{2008}, we conclude that
$F(\phi_0)>0$ on $M$, or else $F(\phi_0)\in B_+$. \ep

{\bf Proof of Theorem \ref{1001}}  \ \  For completeness, let us
redo the argument of \cite{De}. First, Proposition \ref{2007} has
shown that $F|_{A_+}:A_+\rightarrow F(A_+)$ is a diffeomorphism. If
$J$ is integrable, then $(M,g,J,\omega)$ is a closed K\"{a}hler
manifold. The restricted operator $F|_{A_+}$ is a surjectivity map
since there always exists a solution of Calabi-Yau equation on
K\"{a}hler manifold \cite{Ya1}. Suppose that $J$ is not integrable.
From  Lemma \ref{3006}, we have constructed a function $\phi_0$
belonging to the boundary of $A_+$ and $F(\phi_0)>0$ on $M$. We
claim that $F(\phi_0)$ does not belong to the image of $F$ on $A_+$.
Otherwise, we may assume that there exists a function $\phi_1\in
A_+$ such that
$$F(\phi_0)=F(\phi_1).$$
For $t\in (0,1]$, let $$\phi_t=t\phi_1+(1-t)\phi_0.$$ Obviously, for
all $t\in (0,1],$ $\phi_t\in A_+$. Hence
$$\int_0^1\frac{d}{dt}[F(\phi_t)]dt=0.$$
But
$$\int_0^1\frac{d}{dt}[F(\phi_t)]dt\omega^n
=\int_0^1n(\omega(\phi_t))^{n-1}\wedge(d(Jd(\phi_1-\phi_0))dt.$$
 Then, by Lemma \ref{2006}, we have that
 $$L(\phi_0)(f)\omega^n=\int_0^1n(\omega(\phi_t))^{n-1}dt\wedge d(Jdf)$$
is a linear elliptic operator of second order without zero-th term.
Making use of maximal principle of Hopf, we obtain that the kernel
of $L(\phi_0)$ consists of constant functions. Therefore,
$$\phi_0=\phi_1$$ (note that $\phi_0$, $\phi_1$ belong to
$\bar{A}_+$) which implies that
$$\phi_0\in A_+.$$
This contradicts the fact that $\phi_0$ belongs to the boundary of
$A_+$.\\

{\bf Acknowledgements} The authors would like to thank Xifang CAO,
 Jinquan LUO and Ying ZHANG for helpful discussions. The authors
 also express their sincere gratitude toward the referee for having
 pointed out several nice and pertinent improvements in the writing
 of the paper.

\end{document}